\documentclass[10pt]{amsart}

\usepackage{amsmath, amssymb, amsthm, amsfonts, mathrsfs}

\input xy
\xyoption{all}

\newtheorem{prop}{Proposition}[section]
\newtheorem{defn}[prop]{Definition}
\newtheorem{lemma}[prop]{Lemma}
\newtheorem{remark}[prop]{Remark}
\newtheorem{theorem}[prop]{Theorem}
\newtheorem{cor}[prop]{Corollary}
\theoremstyle{definition}

\theoremstyle{definition}

\newcommand{\infinity}{\ensuremath{\infty } }

\begin{document}

\title{On indecomposable trees in the boundary of Outer space}

\author{Patrick Reynolds}

\address{\tt Department of Mathematics, University of Illinois at
  Urbana-Champaign, 1409 West Green Street, Urbana, IL 61801, USA}
  \email{\tt
  preynol3@math.uiuc.edu} 

\date{\today}

\maketitle

\begin{abstract}
Let $T$ be an $\mathbb{R}$-tree, equipped with a very small action of the rank $n$ free group $F_n$, and let $H \leq F_n$ be finitely generated.  We consider the case where the action $F_n \curvearrowright T$ is indecomposable--this is a strong mixing property introduced by Guirardel.  In this case, we show that the action of $H$ on its minimal invarinat subtree $T_H$ has dense orbits if and only if $H$ is finite index in $F_n$.  There is an interesting application to dual algebraic laminations; we show that for $T$ free and indecomposable and for $H \leq F_n$ finitely generated, $H$ carries a leaf of the dual lamination of $T$ if and only if $H$ is finite index in $F_n$.  This generalizes a result of Bestvina-Feighn-Handel regarding stable trees of fully irreducible automorphisms.
\end{abstract}

\section{Introduction}

Let $G$ be a finitely generated group, and suppose that $G \curvearrowright T$ is an action by isometries of $G$ on an $\mathbb{R}$-tree $T$.  

\begin{defn}  \emph{Following \cite{Gui08}, we say that the action $G \curvearrowright T$ is \emph{indecomposable} if for any non-degenerate arcs $I,J \subseteq T$, there are elements $g_1,...,g_r \in G$ such that $J \subseteq g_1I \cup...\cup g_rI$ and such that $g_iI \cap g_{i+1}I$ is non-degenerate for $i \leq r-1$.}  
\end{defn}

It is important to note that the intersections $g_iI \cap g_{i+1}I$ need not be contained in $J$, or even interect $J$ non-degenerately; see \cite{Gui08} for further discussion.  Indecomposability of the action $G \curvearrowright T$ is a strong mixing property; it prohibits the existence of a \emph{transverse family} for the action $G \curvearrowright T$ (see Definition \ref{D.TF}).  In particular, if the action $G \curvearrowright T$ is indecomposable, then $G \curvearrowright T$ cannot be written as a non-trivial \emph{graph of actions} (see \cite{Lev94, Gui08}).  If $H \leq G$ is a finitely generated subgroup containing a hyperbolic isometry of $T$, then there is a canonical minimal subtree $T_H$ for the action $H \curvearrowright T$; notice that if the action $G \curvearrowright T$ has dense orbits, and if $H \leq G$ is a finitely generated, finite index subgroup, then the action $H \curvearrowright T$ has dense orbits as well.  The main result of this paper says that, in some sense, certain indecomposable actions cannot contain any interesting subactions other than the obvious ones. 

Let $\overline{cv}_n$ denote the unprojectivised closed Outer space, \emph{i.e.} the space of very small actions of $F_n$ on $\mathbb{R}$-trees (see Definition~\ref{D.CVn}); we show:

\vspace{.2cm}

\noindent \textbf{Theorem~\ref{T.Main}.}  \emph{Suppose that $T \in \overline{cv}_n$ is indecomposable, and let $H \leq F_n$ be finitely generated.  The action $H \curvearrowright T_H$ has dense orbits if and only if $H$ has finite index in $F_n$.}  

\vspace{.2cm}

There is a nice application of Theorem \ref{T.Main} to \emph{algebraic laminations}: associated to any action $F_n \curvearrowright T$ of $F_n$ on a tree is a \emph{dual lamination} $L^2(T) \subseteq \partial^2(F_n)$, which is an algebraic analog of a surface lamination; here $\partial^2(F_n):=\partial F_n \times \partial F_n - \Delta$ (see section 2.3 for an brief introduction, and \cite{CHL08a, CHL08b} for details).  We say that a finitely generated subgroup $H \leq F_n$ \emph{carries} a leaf $l \in L^2(T)$ if $l \in \partial^2(H) \subseteq \partial^2(F_n)$.  

\vspace{.2cm}

\noindent \textbf{Corollary~\ref{C.Main}.}  \emph{Suppose that $T \in \overline{cv}_n$ is indecomposable and free with dual lamination $L^2(T)$, and let $H \leq F_n$ be finitely generated.  Then $H$ carries a leaf of $L^2(T)$ if and only if $H$ is finite index in $F_n$.}  

\vspace{.2cm}

The reason for the assumption that the action be free comes from the definition of the dual lamination of an action $F_n \curvearrowright T$; namely, if $K \leq F_n$ has a fixed point in $T$, then $\partial^2(K) \subseteq L^2(T)$.  Further, since the action $F_n \curvearrowright T$ is minimal, it is the case that $K$ is infinite index in $F_n$.  

The results of this paper can be thought of as a dynamical-algebraic analogy between indecomposable trees in the boundary of Outer space and ending laminamtions on surfaces.  A lamination $L$ on a compact surface $S$ (possibly with boundary) is \emph{minimal} if every half leaf of $L$ is dense in $L$, and $L$ is \emph{filling} if all complimentary regions are ideal polygons or crowns.  If $L$ is minimal and filling, then $L$ is called an \emph{ending lamination} (see \cite{CB88} for background on suface laminations).  From the definition it is evident that a finite cover of an ending lamination is also an ending lamination.  

In \cite{Sc78} Scott proves that suface groups are subgroup separable (or LERF); his proof is geometric: he finds, for any finitely generated subgroup $H \leq \pi_1(S)$, a finite cover $S_1 \rightarrow S$, a compact surface $S'$, along with a $\pi_1$-injective embedding $\iota: S' \rightarrow S_1$ such that $\pi_1(\iota)(\pi_1(S'))=H$.  This geometric description of subgroups of $\pi_1(S)$ gives a clear picture of which subgroups of $\pi_1(S)$ are able to ``encode'' leaves of the lamination on $S$.  Say that a finitely generated subgroup $H \leq \pi_1(S)$ \emph{carries} a leaf $l$ of $L$ if there are $S_1$, $\iota$, and $S'$ as above, such that a lift of $l$ in $S_1$ is contained in $\iota(S')$.  If $S$ is equipped with an ending lamination $L$, it is evident that the lifted lamination $L_1$ on $S_1$ intersects $\iota(S')$ in finite arcs, unless $\iota(S')=S_1$, \emph{i.e.} unless $H$ is finite index in $\pi_1(S)$.  It follows that no finitely generated subgroup of infinite index carries a leaf of $L$.  Now suppose that $\mathscr{L}=(L, \mu)$ is a measured lamination with $L$ an ending lamination, and let $\mathscr{L}_1=(L_1, \mu_1)$ the lift of $\mathscr{L}$ to $S_1$.  Let $T_{\mathscr{L}}$ denote the $\mathbb{R}$-tree dual to $\mathscr{L}$, and let $T_H \subseteq T$ the minimal invariant subtree for the action of $H$ on $T$.  Evidently, $T_H$ is ``dual'' to $\mathscr{L}_1':=(L_1\cap \iota(S'),\mu_1|_{\iota(S')})$, so the action $H \curvearrowright T_H$ is discrete, again unless $H$ is finite index in $\pi_1(S)$.  Hence, it follows that $H$ carries a leaf of $L$ if and only if $H$ is finite index in $\pi_1(S)$ if and only if the action $H \curvearrowright T_H$ is indiscrete.   

It is easy to see that if an action $\pi_1(S) \curvearrowright T_{\mathscr{L}}$ is dual to a measured ending lamination $\mathscr{L}$ on $S$, then the action is indecomposable.  On the other hand, it follows from Skora's duality theorem \cite{Sko96} and the Rips theory (\cite{BF95}, \cite{GLP94}) that any indecomposable, relatively elliptic action $\pi_1(S) \curvearrowright T$ is dual to an ending lamination on $S$; here \emph{relatively elliptic} means that the (maximal) elliptic subgroups of the action $\pi_1(S) \curvearrowright T$ are precisely the peripheral subgroups of $\pi_1(S)$.  There are other natural examples of indecomposable trees.  The first come from the Rips theory: any geometric tree dual to a minimal band complex is indecomposable (see \cite{BF95} for explanation of terminology and \cite{Gui08} for a proof); this includes the ``surface trees'' mentioned above as well as the so-called thin (or exotic, or Levitt) trees (see \cite{BF95, GLP94, MKap01, Gui08} for details).  Finally, stable trees of fully irreducible (iwip) automorphisms are indecomposable; this can be shown using the machinery of  \cite{BH92} and \cite{BFH97}.  There are examples of such ``iwip trees'' that are not geometric \cite{BF95}.  

As mentioned above, the main results of this paper are known for surface trees.  Using train track machinery, Bestvina-Feighn-Handel establish these results in the special case of stable trees of fully irreducible automorphisms (\cite[Theorem 5.4]{BFH97} and \cite[Proposition 2.4]{BFH97}.  We remark that it follows from the North-South dynamics \cite{LL03} that any stable tree of a fully irreducible automorphism is uniquely ergodic (see Section 3.0 below); on the other hand, \cite{Mar97} establishes the existence of non-uniquely ergodic thin band complexes, so the theorem is saying something new even in the case of geometric trees.  

The inspiration for the proof of the main result is precisely the discussion presented above regarding the dynamical-algebraic properties of ending laminations and their dual trees; in fact, the skeleton of the current proof is essentially identical to that surface theory argument.  The first ingredient is a result of Guirardel, which says that any ``finite cover'' of an indecomposable action $G \curvearrowright T$ is also indecomposable.  We then establish a certain measure-theoretic approximation of actions $F_n \curvearrowright T \in \overline{cv}_n$ with dense orbits: we show that any such action is ``supported almost everywhere'' on a finite forest of arbitrarily small measure, and this allows us to construct from the action $F_n \curvearrowright T$ a finitely generated pseudogroup (see Definition~\ref{D.PG}) with well-controlled dynamics (see Lemma~\ref{L.AESupp}).   All this is combined with an inequality of Gaboriau-Levitt-Paulin to greatly restrict the ``shape'' of families $\{gT_H\}_{g \in F_n}$ for e-algebraically closed subgroups $H \leq F_n$ (see Definition~\ref{D.EAC}).  Finally, the strong subgroup separability of $F_n$ is used to conclude.

The techniques appearing in this paper are part of a more robust approach to studying ``subgroup actions'' in $\overline{cv}_n$; in a forthcoming paper, we refine these techniques to study the dynamics of the action of an irreducible endomorphism of $F_n$ on $\overline{cv}_n$ \cite{R10b}.

\textbf{Acknowledgements}:  Special thanks go to Vincent Guirardel for explaining that finite covers of indecomposable trees are also indecomposable.  Very special thanks go to my advisor Ilya Kapovich for investing a huge amount of time and energy into my graduate education.

The author acknowledges support from National Science Foundation grant DMS 08-38434 ”EMSW21-MCTP: Research Experience for Graduate Students.

\section{Background}

In this section we briefly review the relevant definitions around $\mathbb{R}$-trees, Outer space, and algebraic laminations.  In what follows $F_n$ denotes the free group of rank $n$.

\subsection{Basics About $\mathbb{R}$-Trees}

A metric space $(T,d)$ is called an $\mathbb{R}$-\emph{tree} (or just a \emph{tree}) if for any two points $x,y \in T$, there is a unique topological arc $p_{x,y}:[0,1] \rightarrow T$ connecting $x$ to $y$, and the image of $p_{x,y}$ is isometric to the segement $[0,d(x,y)]$.  As is usual, we let $[x,y]$ stand for Im$(p_{x,y})$, and we call $[x,y]$ the \emph{segment} (also called an \emph{arc}) in $T$ from $x$ to $y$.  A segment is called \emph{non-degenerate} if it contains more than one point.  We let $\overline{T}$ stand for the metric completion of $T$.  Unless otherwise stated, we regard $T$ as a topological space with the metric topology.  If $T$ is a tree, and $x \in T$, then $x$ is called a \emph{branch point} if the cardinality of $\pi_0(T - \{x\})$ is greater than two.  For $x \in T$, the elements of $\pi_0(T- \{x\})$ are called \emph{directions} at $x$.

In this paper, all the trees we consider are equipped with an isometric (left) action of a finitely generated group $G$, i.e. a group morphism $\rho: G \rightarrow$ Isom$(T)$; as usual, we always supress the morphism $\rho$ and identify $G$ with $\rho(G)$.  A tree $T$ equipped with an isometric action will be called an $G$-\emph{tree}, and we denote this situation by $G \curvearrowright T$.  Notice that an action $G \curvearrowright T$ induces an action of $G$ on the set of directions at branch points of $T$.  We identify two $G$-trees $T,T'$ if there is an $G$-equivariant isometry between them.

There are two sorts of isometries of trees: an isometry $g$ of $T$ is called \emph{elliptic} if $g$ fixes some point of $T$, while an isometry $h$ of $T$ is called \emph{hyperbolic} if it is not elliptic.  It is easy to see that any hyperbolic isometry $h$ of $T$ leaves invariant a unique isometric copy of $\mathbb{R}$ in $T$ which is called the \emph{axis} of $h$ and denoted by $A(h)$.  If $g$ is an elliptic isometry, we let $A(g)$ stand for the fixed point set of $g$, \emph{i.e.} $A(h):=\{x \in T|hx=x\}$.  Given a $G$-tree $T$, we have the so-called \emph{hyperbolic length function} $l_T:G \rightarrow \mathbb{R}$, where
$$l_T(g):=\inf \{d(x,gx)|x \in T\}$$ 

\noindent The number $l_T(g)$ is called the \emph{translation length} of $g$, and it is easily verified that, for any $g \in F_N$, the infimum is always realized on $A(g)$, so that $g$ acts on $A(g)$ as a translation of length $l_T(g)$.  If $H \leq G$ is a finitely generated subgroup containing a hyperbolic isometry, then $H$ leaves invariant the set
$$T_H:=\cup_{l_T(h)>0} A(h)$$

\noindent which is a subtree of $T$, and is minimal in the set of $H$-invariant subtrees of $T$; $T_H$ is called the \emph{minimal invariant subtree for} $H$.  An action $G \curvearrowright T$ is called \emph{minimal} if $T=T_G$; a minimal action $G \curvearrowright T$ is \emph{non-trivial} if $T$ contains more than one point.

\subsection{Outer Space(s)}

Recall that an action $F_n \curvearrowright T$ is \emph{free} if for any $1 \neq g \in F_n$ one has $l_T(g) > 0$.  If $X \subseteq T$, then the \emph{stabilizer} of $X$ is $Stab(X):=\{g \in F_n|gX=X\}$--the setwise stabilizer of $X$.  We say that an action $F_n \curvearrowright T$ is \emph{very small} if:

\begin{enumerate}
 \item [(i)] $F_n \curvearrowright T$ is minimal,
 \item [(ii)] for any non-degenerate arc $I \subseteq T$, $Stab(I) = \{1\}$ or $Stab(I)$ is a maximal cyclic subgroup of $F_n$,
 \item [(iii)] stabilizers of tripods are trivial.
\end{enumerate}
 
An action $F_n \curvearrowright T$ is called \emph{discrete} (or \emph{simplicial}) if the $F_n$-orbit of any point of $T$ is a discrete subset of $T$; in this case $T$ is obtained by equivariantly assigning a metric to the edges of a (genuine) simplicial tree.  It is important to note that the metric topology is weaker than the simplicial topology if the tree is not locally compact.

Let $T,T'$ be trees; a map $f:T \rightarrow T'$ is called a \emph{homothety} if $f$ is $F_n$-equivariant and bijective, and if there is some positive real number $\lambda$ such that for any $x,y \in T$, we have $d_{T'}(f(x),f(y))=\lambda d_T(x,y)$; in this case $T,T'$ are called \emph{projectively equivalent} or \emph{homothetic}.  

\begin{defn}\label{D.CVn}
\noindent
\begin{enumerate}
 \item The \emph{unprojectivised Outer space} of rank $n$, denoted $cv_n$, is the topological space whose underlying set consists free, minimal, discrete, isometric actions of $F_n$ on $\mathbb{R}$-trees; it is equipped with the \emph{length function topology}.
 \item \cite{CV86} The \emph{Culler-Vogtmann} \emph{Outer space} of rank $n$, denoted $CV_n$, is the topological space whose underlying set consists of homothety classes of free, minimal, discrete, isometric actions of $F_n$ on $\mathbb{R}$-trees; it is equipped with the \emph{projective length function topology}.
 \item The \emph{unprojectivised closed Outer space} of rank $n$, denoted $\overline{cv}_n$, is the topological space whose underlying set consists of very small isometric actions of $F_n$ on $\mathbb{R}$-trees; it is equipped with the \emph{length function topology}.
 \item The \emph{closed Outer space} of rank $n$, denoted $\overline{CV}_n$, is the topological space whose underlying set consists of homothety classes of very small isometric actions of $F_n$ on $\mathbb{R}$-trees; it is equipped with the \emph{projective length function topology}.
\end{enumerate}

\end{defn}

As it is well-known that a minimal $F_n$-tree is completely determined by its hyperbolic length function (see, for example, \cite{Chi}), points in $CV_n$ can be thought of as projective classes of such length functions, \emph{i.e.} $CV_n \subseteq \mathbb{P}\mathbb{R}^{F_n}$; and $CV_n$ is topologized via the quotient of the weak topology on length functions.  It is the case that the closure $\overline{CV}_n$ of $CV_n$ is compact and consists precisely of homothety classes of very small $F_N$-actions on $\mathbb{R}$-trees \cite{CL95, BF94}.  For more background on $CV_n$ and its closure, see \cite{Vog02} and the references therein.

\subsection{Dual Laminations}

Here, we present a brief and restricted view of dual laminations of $F_n$-trees; see \cite{CHL08a} and \cite{CHL08b} for a careful development of the general theory.  Let $\partial F_n$ denote the Gromov boundary of $F_n$--\emph{i.e.} the Gromov boundary of any Cayley graph of $F_n$; let $\partial^2(F_n):=\partial F_n \times \partial F_n - \Delta$, where $\Delta$ is the diagonal.  The left action of $F_n$ on a Cayley graph induces actions by homeomorphisms of $F_n$ on $\partial F_n$ and $\partial^2 F_n$.  Let $i: \partial^2 F_n \rightarrow \partial^2 F_n$ denote the involution that exchanges the factors.  An \emph{algebraic lamination} is a non-empty, closed, $F_n$-invariant, $i$-invariant subset $\mathscr{L} \subseteq \partial^2 F_n$.  

Fix an action $F_n \curvearrowright T$ with dense orbits; following \cite{LL03} (see also \cite{CHL08b}), we associate an algebraic lamination $L^2(T)$ to the action $F_n \curvearrowright T$. Let $T_0 \in cv_n$ (\emph{i.e.} the action $F_n \curvearrowright T_0$ is free and discrete), and let $f:T_0 \rightarrow T$ be an $F_n$-equivariant map, isometric when restricted to edges of $T_0$.  Say that $f$ has \emph{bounded backtracking} if there is $C > 0$ such that $f([x,y]) \subseteq N_C([f(x),f(y)])$, where $N_C$ denotes the $C$-neighborhood.  For $T_0 \in cv_n$, denote by $vol(T_0):=vol(T_0/F_n)$ the sum of lengths of edges of the finite metric graph $T_0/F_n$.

\begin{prop}\label{P.BBT}\cite[Lemma 2.1]{LL03}
Let $T \in \overline{cv}_n$; let $T_0 \in cv_n$; and let $f:T_0 \rightarrow T$ be equivariant and isometric on edges.  Then $f$ has bounded backtracking with $C=vol(T_0)$.  
\end{prop}
 
For $T_0 \in cv_n$, we have an identification $\partial T_0 \cong \partial F_n$.  If $\rho$ is a ray in $T_0$ representing $X \in \partial F_n$, we say that $X$ is $T$-\emph{bounded} if $f \circ \rho$ has bounded image in $T$; this does not depend on the choice of $T_0$ (see \cite{BFH97}).  

\begin{prop}\label{P.DefQ}\cite[Proposition 3.1]{LL03}
Let $T \in \overline{cv}_n$ have dense orbits, and suppose that $X \in \partial F_n$ is $T$-bounded.  There there is a unique point $Q(X) \in \overline{T}$ such that for any $f:T_0 \rightarrow T$, equivariant and isometric on edges, and any ray $\rho$ in $T_0$ representing $X$, the point $Q(X)$ belongs to the closure of the image of $f \circ \rho$ in $\overline{T}$.  Further the image of $f \circ \rho$ is a bounded subset of $\overline{T}$.  
\end{prop}

The (partially-defined) map $Q$ given above is clearly $F_n$-equivariant; in fact, it extends to an equivariant map $Q: \partial F_n \rightarrow \overline{T} \cup \partial T$, which is surjective (see \cite{LL03}).  The crucial property for us is that $Q$ can be used to associate to $T$ an algebraic lamination.  

\begin{prop}\label{P.DefL2}\cite{CHL08b}
Let $T \in \overline{cv}_n$ have dense orbits.  The set $L^2_Q(T):=\{(X,Y) \in \partial^2(F_n)|Q(X)=Q(Y)\}$ is an algebraic lamination.
\end{prop}

Following \cite{CHL08b}, we mention that there is different, perhaps more intuitive, procedure for defining $L^2(T)$.  Let $T \in \overline{cv}_n$ (not necessarily with dense orbits, but not free and discrete), and let $\Omega_{\epsilon}(T):=\{g \in F_n|l_T(g) < \epsilon\}$, where $l_T$ is the hyperbolic length function for the action $F_n \curvearrowright T$.  The set $\Omega_{\epsilon}(T)$ generates an algebraic lamination $L^2_{\epsilon}(T)$, which is the smallest algebraic lamination containing $(g^{-\infinity},g^{\infinity})=(...g^{-1}g^{-1},gg...) \in \partial^2(F_n)$ for every $g \in \Omega_{\epsilon}$.  One then defines $L^2_{\Omega}(T):= \cap_{\epsilon > 0} L^2_{\epsilon}(T)$.  In \cite{CHL08b} it is shown that for an action $F_n \curvearrowright T \in \overline{cv}_n$ with dense orbits, $L^2_{\Omega}(T)=L^2_Q(T)$, as defined above.  

\begin{defn}\label{D.L2}
Let $F_n \curvearrowright T \in \overline{cv}_n$ be an action with dense orbits.  The \emph{dual lamination of} $F_n \curvearrowright T$ is $L^2(T):=L^2_Q(T)=L^2_{\Omega}(T)$.  
\end{defn}

\section{Invariant Measures and Transverse Families}

Let $T$ be an $\mathbb{R}$-tree.  

\begin{defn}\label{D.LM}\cite{Gui00}
A \emph{length measure} (or just \emph{measure}) $\mu$ on $T$ is a collection $\mu=\{\mu_I\}_{I\subseteq T}$ of finite positive Borel measures on the finite arcs $I \subseteq T$; it is required that for $J \subseteq I$, $\mu_J=(\mu_I)|_J$.   
\end{defn}

As these measures are defined locally on finite arcs, all the usual measure-theoretic definitions are similarly defined: a set $X \subseteq T$ is $\mu$-\emph{measurable} if $X \cap I$ is $\mu_I$-measurable for each $I \subseteq T$; $X$ has $\mu$-\emph{measure zero} if $X \cap I$ is $\mu_I$-measure zero for each $I$; and so on.  The \emph{Lebesgue length measure}, denoted $\mu_L$, on $T$ is the collection of Lebesgue measures on the finite arcs of $T$.  If $T$ is equipped with an action of a group $G$, then we say that a measure $\mu$ is $G$-\emph{invariant} if $\mu_I(X \cap I) = \mu_{gI}(gX \cap gI)$ holds for each $\mu$-measurable set $X$ and each $g \in G$.  Note that if the action $G \curvearrowright T$ is by isometries, then the Lebesgue measure is invariant.  We let $M(T)=M(G \curvearrowright T)$ stand for the set of invariant measures on $T$.

Suppose that $G \curvearrowright T$ is an action by isometries, with $G$ a countable group.  Say that the action is \emph{finitely supported} if there is a finite subtree $K \subseteq T$ such that any finite arc $I \subseteq T$ may be covered by finitely many translates of $K$ by elements of $G$; in this case, we say that the action $G \curvearrowright T$ is \emph{supported on} $K$.  Note that, if $G$ is finitely generated, then any minimal action $G \curvearrowright T$ is finitely supported.  

Let $K$ be a compact topological space. 

\begin{defn}\label{D.PG}
A collection of partially defined homeomorphisms $\Gamma$ of $K$ is called a \emph{pseudogroup} if the following are satisfied:

\begin{enumerate}
 \item [(1)] the identity mapping is an element of $\Gamma$,
 \item [(2)] if $\gamma \in \Gamma$, then $\gamma^{-1} \in \Gamma$, where $Dom(\gamma^{-1})=Ran(\gamma)$,
 \item [(3)] if $\gamma_1, \gamma_2 \in \Gamma$, then $\gamma_1 \circ \gamma_2 \in \Gamma$
 \item [(4)] if $\gamma_1, \gamma_2 \in \Gamma$, and if $\gamma_1(x)=\gamma_2(x)$ for all $x \in Dom(\gamma_1) \cap Dom(\gamma_2)$, then if $\gamma_1 \cup \gamma_2$ is a homeomorphism, then $\gamma_1 \cup \gamma_2 \in \Gamma$, and
 \item [(5)] if $\gamma_1 \in \Gamma$, then the restriction of $\gamma_1$ to any Borel subset of $Dom(\gamma_1)$ is in $\Gamma$.
\end{enumerate}
\end{defn}

We say that $\{\gamma_1,...,\gamma_k,...\}$ \emph{generate} $\Gamma$ if any $\gamma \in \Gamma$ can be obtained from the $\gamma_i$ via the operations in the definition of a pseudogroup.  A measure $\mu$ on $K$ is said to be $\Gamma$-\emph{invariant} if for any measurable $X \subseteq K$, we have $\mu(X \cap dom(\gamma)) = \mu(\gamma (X \cap dom(\gamma)))$ for each $\gamma \in \Gamma$.  We let $M(K)=M(\Gamma,K)$ stand for the set of invariant measures on $K$.  

Let $G \curvearrowright T$ be an action supported on the finite subtree $K \subseteq T$.  We consider the (countably generated) pseudogroup $\Gamma:=\{g|_{K'}:g \in G, K' \subseteq K\}$ of restrictions of the isometries $G$ to Borel subsets of $K$.  Since the action is supported on $K$, there is a bijective correspondence between $M(T)$ and $M(K)$.

A non-trivial measure $\mu \in M(T)$ is called \emph{ergodic} if any $G$-invariant subset is either full measure or zero measure.  A $G$-tree $T$ is called \emph{uniquely ergodic} if there is a unique, up to scaling, $G$-invariant measure $\mu$ on $T$; in this case $\mu$ must be ergodic.  Let $M_0(T)$ denote the set of non-atomic, $G$-invariant measures on $T$, and let $M_1(T) :=\{\nu \in M_0(T)|\nu \leq \mu_L\}$.  Note that both $M_0(T)$ and $M_1(T)$ are convex.  

\begin{prop}\label{P.M0T}\cite[Corollary 5.4]{Gui00}
Let $T \in \overline{cv}_n$ be with dense orbits. Then $M_0(T)$ is a finite dimensional convex set, which is projectively compact. Moreover, $T$ has at most $3n-4$ non-atomic ergodic measures (up to homothety), and every measure in $M_0(T)$ is a sum of these ergodic measures.  Further $M_1(T)$ is compact.
\end{prop}

\subsection{Finite Systems of Isometries}

A \emph{finite tree} is a tree that is the convex hull of a finite set; a \emph{finite forest} is a finite union of finite trees.  A \emph{finite pseudogroup} is a finitely generated pseudogroup $S=(F,A)$, where $F$ is a finite forest.  Let $S=(F,A)$ be a finite pseudogroup generated by $A=\{a_1,...,a_n\}$; we require that $dom(a_i)$ be closed.  For $a_i \in A$, let $B_i:=dom(a_i) \times I$; regard $B_i$ as foliated by leaves of the form $\{pt.\} \times I$.  Form the \emph{suspension} $\Sigma(S)$ of $S$ from the disjoint union $K \sqcup B_1 \sqcup ... \sqcup B_n$ by identifying $B_i \times \{0\}$ with $dom(a_i)$ and $y=(x,1) \in B_i \times \{1\}$ with $a_i(x)$.  Put a relation $R_l$ on points of $\Sigma(S)$, where $x, y \in R_l$ if and only if $x,y$ are contained in a leaf of some $B_i$; let $\overline{R_l}$ be the smallest equivalence relation containing $R_l$; and regard $\Sigma(S)$ as foliated by leaves that are the classes of $\overline{R_l}$.  Note that for $x \in K$, the leaf $l(x)$ containing $x$ intersects $K$ precisely in the orbit $S.x$.  

Let $B$ denote the set of branch points of $K$, and let $E$ denote the set containing all endpoints of all $dom(a_i)$; put $C:=B \cup E$.  A leaf $l$ of $\Sigma(S)$ is called \emph{singular} if $l \cap C \neq \emptyset$; any leaf that is not singular is called \emph{regular}.  Suppose that $\Sigma(S)$ contains a finite regular leaf $l=l(x)$, then for $y \in K$ close to $x$, $l(y)$ is finite and regular.  It follows that there are $y_1,y_2 \in K$ with $x \in [y_1,y_2]$ and $d(y_1,y_2)$ maximal, such that for $z \in (y_1,y_2)$, $l(z)$ is finite and regular.  Hence, $F_x:=\cup_{z \in (y_1,y_2)}l(z)$ is a $(y_1,y_2)$-bundle over some leaf $l(z) \in F_x$.  The set $F_x$ is called a \emph{maximal family of finite orbits}, and the \emph{transverse measure} of $F_x$ is $d(y_1,y_2)$.  Evidently, $l(y_i)$ are singular, so there are finitely many maximal families of finite orbits in $\Sigma(S)$.  This gives a coarse decomposition of $\Sigma(S)$, which is the starting point for a refined decomposition of $\Sigma(S)$, see \cite{GLP94} for the statement as well as for details regarding the above discussion.  

Suppose that $S=(F,A)$ is a finite pseudogroup; define the following:
\begin{enumerate}
 \item $m:=$ total measure of $F$
 \item $d:=$ the sum of measures of domains of generators
 \item $e:=$ the sum of transverse measures of maximal families of finite orbits.  
\end{enumerate}

We regard $m$, $d$, and $e$ as functions $\{$finte pseudogroups$\} \rightarrow \mathbb{R}$.  Say that $S$ has \emph{independent generators} if no reduced word in the generators $A$ and their inverses defines a partial isometry of $F$ that fixes a non-degenerate arc.  

\begin{prop}\label{P.GLP}\cite[Proposition 6.1]{GLP94}
Let $S$, $F$, $A$ as above, and suppose that $S$ has independent generators, then $e(S)+d(S)=m(S)$.
\end{prop}

\subsection{Transverse Families}

Fix a basis $A$ for $F_n$ and an action $F_n \curvearrowright T \in \overline{cv}_n$ with dense orbits, and let $\mu \in M(T)$.  

\begin{defn}\label{D.AESupp}
Say that the action $F_n \curvearrowright T$ is \emph{supported} $\mu$-\emph{a.e.} on a $\mu$-measurable set $X \subseteq T$ if for any arc $I \subseteq T$ and any $\delta > 0$, there are $g_1,...,g_r \in F_n$ such that $\mu(I - (g_1X \cup ... \cup g_rX)) < \delta$.  
\end{defn}

For a finite forest $F \subseteq T$, we write $S=(F,A)$ to denote the pseudogroup generated by restrictions of elements of $A$ to $F$.  Recall that $\mu_L$ denotes Lebesgue measure on $T$.

\begin{lemma}\label{L.AESupp}
Let $T \in \overline{cv}_n$ be with dense orbits.  For any $\epsilon > 0$ and any finite forest $K \subseteq T$, there are finite forests $F_{\epsilon}$ and $F$ such that:
\begin{enumerate}
 \item [(i)] $\mu_L(F_{\epsilon}) < \epsilon$,
 \item [(ii)] $F_n \curvearrowright T$ is supported $\mu_L$-a.e. on $F_{\epsilon}$
 \item [(iii)] $\mu_L(F \cap K) > \mu_L(K) - \epsilon$,
 \item [(iv)] $S=(F,A)$ satisfies $m(S) - d(S) < \epsilon$
\end{enumerate}
\end{lemma}

\begin{proof}
Let $T$, $\epsilon$, and $K$ as in the statement.  By Proposition~\ref{P.M0T} we have that $\mu_L= \sum_{i=1}^p \nu_i$, with each $\nu_i$ ergodic.  Take $J_i \subseteq T$ finite arcs such that $\nu_i(J_i) > 0$ and $\sum_i \mu_L(J_i) < \epsilon$; put $F_{\epsilon}:= \cup_i J_i$.  By ergodicity of the measures $\nu_i$, we get that $\bigcup_{g \in F_n} gF_{\epsilon}$ is a full measure subset of $T$, so the action $F_n \curvearrowright T$ is supported $\mu_L$-a.e. on $F_{\epsilon}$.  Hence, there are $g_1,...,g_r \in F_n$ such that $\mu_L(K \cap (\cup_i g_iF_{\epsilon})) > \mu_L(K) - \epsilon$.  Let $S_0:=(F_{\epsilon},A)=(F_0,A)$ the finite pseudogroup of restrictions of elements of $A$ to $F_0$, and note that $m(S_0) -d(S_0) \leq \mu_L(F_{\epsilon}) < \epsilon$.  Define $F_i:=F_{i-1} \cup \bigcup_{a \in A^{\pm}} aF_{i-1}$, and let $S_i=(F_i,A)$.  Immediately, one has $m(S_i) - d(S_i) \leq m(S_{i-1}) - d(S_{i-1})$.  The claim follows by observing that for $i \geq max\{|g_i|_A\}$ we have $\cup_i g_iF_0 \subseteq F_i$.
\end{proof}

\begin{remark}\label{R.AESupp}
If an action $F_n \curvearrowright T \in \overline{cv}_n$ is indiscrete, but not with dense orbits, then $T$ splits as a graph of actions with vertex trees either finite arcs of $T$ (the \emph{simplicial part} of $T$) or subtrees $T_v$ such that the action $Stab(T_v) \curvearrowright T_v$ is with dense orbits (see \cite{Lev94, Gui08}).  In this case, it follows from the above argument that for any $\epsilon > 0$, the action $F_n \curvearrowright T$ is supported $\mu_L$-a.e. on a finite forest $F_{\epsilon}'$ with $\mu_L(F_{\epsilon}') < vol(T/F_n) + \epsilon$, where $vol(T/F_n)=\inf \mu_L(S)$ with the infimum taken over all measurable $S \subseteq T$ projecting onto $T/F_n$ under the natural map.
\end{remark}  

\begin{defn}\label{D.EAC}\cite{MVW07}
A finitely generated subgroup $H \leq F_n$ is \emph{e-algebraically closed} if for any $g \in F_n - H$, one has $\langle H,g \rangle \cong H \ast \langle g \rangle$.
\end{defn}

Equivalently, $H$ is e-algebraically closed if there is no non-trivial equation $w(\overline{h},x)$ over $H$ with a solution $w(\overline{h},g)$ for $g \in F_n - H$.  Any free factor of $F_n$ is necessarily e-algebraically closed; further, if $H \leq F_n$ has rank $r$ and is maximal in the poset of rank $r$ subgroups of $F_n$, then $H$ is e-algebraically closed (see \cite{MVW07} for details).  

\begin{defn}\label{D.TF}
Let $T$ be an $\mathbb{R}$-tree, equipped with an action of a group $G$.  A $G$-invariant collection $\{T_v\}_{v \in V}$ of non-degenerate proper subtrees of $T$ is called a \emph{transverse family} if whenever $T_v \neq T_{v'}$, $T_v \cap T_{v'}$ contains at most one point.  
\end{defn}

\begin{lemma}\label{L.TF}
Let $T \in \overline{cv}_n$ be with dense orbits; let $H \leq F_n$ a finitely generated subgroup with minimal invariant tree $T_H \subseteq T$.  Suppose that the action $H \curvearrowright T_H$ has dense orbits and that $H$ is e-algebraically closed.  The family of translates $\{gT_H\}_{g \in F}$ is a transverse family.
\end{lemma}

\begin{proof}
Let $T$ and $H$ as in the statement of the lemma.  Note that since $H$ is e-algebraically closed, if $H$ is a proper subgroup of $F_n$, then $H$ is infinite index in $F_n$.  If $H=F_n$, then the statement is trivial, so we suppose that $H$ has infinite index in $F_n$.  Choose a basis $B$ for $H$.  Let $F \subseteq T_H$ be a finite forest; since the action $H \curvearrowright T_H$ is very small, it is the case that the pseudogroup $S=(F,B)$ generated by restrictions of the elements of $B$ to $F$ has independent generators; further since $H$ is e-algebraically closed, it is the case that for any $f \in F_n - H$, the restrictions of $B \cup \{f\}$ to $F$ give a finite pseudogroup with independent generators.  

Toward a contradiction, suppose that there is $f \in F - H$ such that $fT_H \cap T_H$ contains more than one point.  Since the intersection of two trees is convex, we have that $fT_H \cap T_H$ contains a non-degenerate arc $I$.  Choose $\epsilon > 0$ small with respect to $\mu_L(I)$.  Set $K:=I \cup f^{-1}I$; by Lemma~\ref{L.AESupp}, we may find a finite forest $F \subseteq T_H$ such that $\mu_L(F \cap K) > \mu_L(K) - \epsilon$ and such that $S=(K,B)$ satisfies $m(S)-d(S) < \epsilon$.  

Now, consider $S':=(K,B \cup \{f\})$; as noted above, $S'$ has independent generators.  On the other hand, it is clear from the construction that $m(S')-d(S') < 0$, a contradiction to Proposition~\ref{P.GLP}.  It follows that $fT_H \cap T_H$ contains at most one point for each $f \in F_n - H$, so $\{gT_H\}_{g \in F_n}$ is a transverse family.
\end{proof}

\begin{remark}\label{R.TF}
The proof actually shows something stronger: if the action $H \curvearrowright T_H$ is indiscrete and if $H$ is e-algebraically closed, then for $g \in F_n$ with $gT_H \neq T_H$, we have that any non-degenerate intersection $gT_H \cap T_H$ is contained in the simplicial part of $T_H$ (See Remark~\ref{R.AESupp}).

Further, the proof shows that if $H$ is e-algebraically closed and if $H$ is a proper subgroup of $F_n$, then $T_H$ is a proper subtree of $T$; later (see Remark~\ref{R.Proper}), we will see that for $H$ finitely generated, $T_H=T$ if and only if $H$ is finite index in $F_n$.
\end{remark}

\section{Indecomposable Trees}

Recall that a $G$-tree $T$ is called \emph{indecomposable} if for any non-degenerate arcs $I,J \subseteq T$, there are elements $g_1,...,g_r$ such that $J \subseteq g_1I \cup ... \cup g_rI$, and $g_iI \cap g_{i+1}I$ is non-degenerate for $i \leq r-1$.    

\begin{lemma}\label{L.TFInd}
If an action $G \curvearrowright T$ is indecomposable, then there is no transverse family for the action $G \curvearrowright T$.
\end{lemma}

\begin{proof}
Suppose that the action $G \curvearrowright T$ is indecomposable; and, toward contradiction, suppose that $\{T_v\}_{v \in V}$ is a transverse family for the action $G \curvearrowright T$.  Recall that each $T_v$ is a proper, non-degenerate subtree of $T$ and that the collection $\{T_v\}_{v \in V}$ is $G$-invariant.  Hence, we may find distinct $T_v$, $T_{v'}$ along with an arc $I \subseteq T$ such that $I \cap T_v$ and $I \cap T_{v'}$ are non-degenerate.  Define $I_0 := I \cap T_v$; by indecomposability of the action $G \curvearrowright T$, there are $g_1,...,g_r \in G$ such that $I \subseteq g_1I_0 \cup ... \cup g_rI_0$ with $g_iI_0 \cap g_{i+1}I_0$ non-degenerate.  Since $\{T_v\}_{v \in V}$ is a transverse family, it follows that $g_iI_0 \subseteq T_v$ for each $i$, hence $T_v=T_{v'}$, a contradiction.  
\end{proof}

\subsection{Lifting Indecomposability}

The idea for the proof of the following lemma was communicated to us by Vincent Guirardel.

\begin{lemma}\label{L.LiftInd}[Guirardel]
Suppose that the action $G \curvearrowright T$ is indecomposable and that $H \leq G$ is a finitely generated and finite index.  Then the action $H \curvearrowright T$ is indecomposable.
\end{lemma}

\begin{proof}
We remark that since $H \leq G$ is finite index, $T_H=T$; without loss, we may assume that $H$ is normal.  For an arc $I \subseteq T$, define a subtree $Y_I \subseteq T$ as follows.  Put $Y_0:=I$, and define $Y_{i+1}:= Y_i \cup \bigcup_h hI$, where the union is taken over elements $h \in H$ such that $gI \cap Y_i$ is non-degenerate.  Finally set $Y_I:= \cup_i Y_i$.  Toward a contradiction assume that the action $H \curvearrowright T$ is not indecomposable; it follows that we may find a non-degenerate arc $I \subseteq T$ such that $Y_I \subsetneq T$.  By construction, the collection $\{hY_I\}_{h \in H}$ is a transverse family for the action $H \curvearrowright T$.   

Let $\{1=g_1,...,g_l\}$ be a left transversal to $H$ in $G$, and let $[g_i]$ denote the coset corresponding to $g_i$.  Consider the collections $\mathscr{Y}_i:=\{gY_I|g \in [g_i]\}$; we claim that there are $Y_i \in \mathscr{Y}_i$ such that $\cap_i Y_i$ is non-degenerate.  Note that by indecomposibility of the action $G \curvearrowright T$, there is $g \in G - H$ such that $gY_I \cap Y_I$ is non-degenerate and $gY_I \neq Y_I$; say $g \in [g_i]$.  Consider the collection of non-degenerate intersections $hY_I \cap gY_I$ for $g \in [g_i]$ and $h \in H$.  This collection is a transverse family for the action $H \curvearrowright T$; indeed, normality of $H$ ensures invariance, so suppose that $(g_ih_1Y_I \cap h_2Y_I) \cap h(g_ih_1Y_I \cap h_2Y_I)$ is non-degenerate.  We have:
\begin{align*}
h(g_ih_1Y_I \cap h_2Y_I) \cap (g_ih_1Y_I \cap h_2Y_I) &= g_ih'h_1Y_I \cap hh_2Y_I \cap g_ih_1Y_I \cap h_2Y_I\\
&=g_i(h'h_1Y_I \cap h_1Y_I) \cap (hh_2Y_I \cap h_2Y_I)
\end{align*}

As $\{hY_I\}_{h \in H}$ is a transverse family for the action $H \curvearrowright T$, it follows that $h'h_1Y_I=h_1Y_I$ and $hh_2Y_I = h_2Y_I$.  Hence $h(g_ih_1Y_I \cap h_2Y_I)=(g_ih_1Y_I \cap h_2Y_I)$.

Again by indecomposability of the action $G \curvearrowright T$, that there are $g' \in G$ and $g \in [g_i]$, with $g'(hY_I \cap gY_I) \cap (Y_I \cap gY_I)$ non-degenerate and $g'(hY_I \cap gY_I) \neq (hY_I \cap gY_I)$.  Evidently, $g \notin [g_i]$, hence we may continue in this way to get a non-degenerate intersection $g_1'Y_I \cap ... \cap g_l'Y_I$, for $g_i' \in [g_i]$.  As above this is a transverse family for the action $H \curvearrowright T$; we claim that it is a transverse family for the action $G \curvearrowright T$.  We have $Y=g_1h_1Y_I \cap ... \cap g_lh_lY_I$ non-degenerate.  Let $g \in G$, then by normality of $H$, we get
\begin{align*}
gY \cap Y &=g(g_1h_1Y_I \cap ... \cap g_lh_lY_I) \cap (g_1h_1Y_I \cap ... \cap g_lh_lY_I)\\
&=g_ih(g_1h_1Y_I \cap ... \cap g_lh_lY_I) \cap (g_1h_1Y_I \cap ... \cap g_lh_lY_I)\\
&=(g_ig_1h_1'h_1Y_I \cap ... \cap g_ig_lh_l'h_lY_I) \cap (g_1h_1Y_I \cap ... \cap g_lh_lY_I)\\
&=(g_{i_1}h_1'h_1Y_I \cap ... \cap g_{i_l}h_l'h_lY_I) \cap (g_1h_1Y_I \cap ... \cap g_lh_lY_I)
\end{align*}

Hence $\{gY\}_{g \in G}$ is a transverse family, a contradiction to Lemma~\ref{L.TFInd}.  
\end{proof}

\subsection{Proof of the Main Result}

The following strong separability result is of central importance to us; for a particularly beautiful proof, see \cite{St83} (see \cite{KM02} for extensions of the ideas of \cite{St83}).  

\begin{theorem}\label{T.MH}(Marshall Hall's Theorem)\cite{MH50}
Let $H \leq F_n$ be finitely generated, and let $g \in F_n - H$.  There is finitely generated $F' \leq F_n$ of finite index, such that $F'=H \ast K$ with $g \notin F'$.
\end{theorem}

Recall that if $F_0 \leq F_n$ is a free factor, then $F_0$ is e-algebraically closed in $F_n$.  In light of this, the above theorem states that for any finitely generated $H \leq F_n$, we can find a finitely generated, finite index subgroup $F' \leq F$ such that $H$ is e-algebraically closed in $F'$ (here, we do not use the subgroup separability).

\begin{theorem}\label{T.Main}
Let $T \in \overline{cv}_n$ be indecomposable, and let $H \leq F_n$ be a finitely generated subgroup.  The action $H \curvearrowright T_H$ is indiscrete if and only if $H$ is finite index in $F_n$.  
\end{theorem}

\begin{proof}
Let $T \in \overline{cv}_n$ be indecomposable, and let $H \leq F_n$ be finitely generated and infinite index in $F_n$.  Toward a contradiction suppose that the action $H \curvearrowright T_H$ is indiscrete.  It follows that there is finitely generated $H' \leq H$ such that the action $H' \curvearrowright T_{H'}$ is with dense orbits, so we may suppose that the action $H \curvearrowright T_H$ is with dense orbits.  By Theorem~\ref{T.MH}, there is a finitely generated $F' \leq F_n$ of finite index, such that $H \leq F'$ is a free factor; hence $H$ is e-algebraically closed in $F'$.  By Lemma~\ref{L.LiftInd} we have that the action $F' \curvearrowright T$ is indecomposable; on the other hand, by Lemma~\ref{L.TF} as $H$ is e-algebraically closed in $F'$, the family of $F'$-translates of $T_H$ is a transverse family for the action $F' \curvearrowright T$.  By Lemma~\ref{L.TFInd}, we arrive at a contradiction to the indecomposability of the action $F' \curvearrowright T$.
\end{proof}

\begin{remark}\label{R.Proper}
A similar line of reasoning as above shows that for any finitely generated, infinite index $H \leq F_n$, we have that $T_H$ is a proper subtree of $T$.  Indeed, by Theorem~\ref{T.MH} we may find $F' \leq F_n$, finite index, such that $H$ is e-algebraically closed in $F'$.  It follows from Remark~\ref{R.TF} that $T_H \subsetneq T_{F'}=T$.  
\end{remark}

\begin{cor}
Let $F_n \curvearrowright T \in \overline{cv}_n$ be an action with dense orbits, and let $H \leq F_n$ be finitely generated.  Then $T_H=T$ if and only if $H$ is finite index in $F_n$.
\end{cor}

To complete the analogy with the dynamical-algebraic properties of ending laminations, we bring Corollary~\ref{C.Main} below; as mentioned in the introduction, the hypothesis that the action $F_n \curvearrowright T$ be free is essential--it is a by-product of the definition of the dual lamination of a tree.  Corollary~\ref{C.Main} follows immediately from Theorem~\ref{T.Main} and the following:

\begin{lemma}\label{L.Carry}
Let $T \in \overline{cv}_n$ be free with dense orbits, and let $H \leq F_n$ finitely generated.  The action $H \curvearrowright T_H$ is indiscrete if and only if $H$ carries a leaf of $L^2(T)$.  
\end{lemma}

\begin{proof}
Suppose that the action $H \curvearrowright T_H$ is indiscrete, then $L^2(T_H):=L^2(H \curvearrowright T_H)$ is non-empty; from the definition of $L^2(T)$, it is evident $L^2(T_H) \subseteq L^2(T)$.

Conversely, suppose that $H$ carries a leaf $l \in L^2(T)$.  Toward a contradiction suppose that the action $H \curvearrowright T_H$ is discrete.  Let $T^0 \in cv_n$, and choose an $F_n$-equivariant map $f: T^0 \rightarrow T$; then $f$ restricts to an $H$-equivarant map $f_H :T^0_H \rightarrow T_H$, which descends to $\overline{f}_H:T^0_H/H \rightarrow T_H/H$, which is a homotopy equivalence, since the action $H \curvearrowright T_H$ is free.  It follows that $f_H$ is a quasi-isometry.  On the other hand, by Proposition~\ref{P.DefQ}, if $l \in L^2(T)$ is carried by $H$, there is a line $l_0 \subseteq T^0_H$ representing $l$ that is mapped via $f_H$ to a bounded subset of $T$, a contradiction.
\end{proof}

\begin{cor}\label{C.Main}
Suppose that $T \in \overline{cv}_n$ is indecomposable and free with dual lamination $L^2(T)$, and let $H \leq F_n$ be finitely generated.  Then $H$ carries a leaf of $L^2(T)$ if and only if $H$ is finite index in $F_n$.  
\end{cor}

\bibliographystyle{amsplain}
\bibliography{indecompREF.bib}
\end{document}